\documentclass[10pt]{article}
\usepackage{amscd,amsmath,latexsym}
\usepackage{latexsym}
\usepackage{amsfonts}
\usepackage{amssymb,amsthm}
\usepackage{longtable}
\usepackage{epsfig}

\theoremstyle{definition}
\newtheorem{Def}{Definition}[section]

\newtheorem{thm}[Def]{Theorem}
\newtheorem{cor}[Def]{Corollary}
\newtheorem{lem}[Def]{Lemma}

\newtheorem{rem}[Def]{Remark}

\newtheorem{p}[Def]{Problem}

\newcommand{\constant}{{\rm constant}}

\newcommand{\re}{{\rm Re}}

\begin{document}
\title{Minimal annuli with constant contact angle\\ along the planar boundaries}
\author{Juncheol Pyo}

\date{}
\maketitle

\renewcommand{\thefootnote}{}
\smallskip
\begin{abstract}
We show that an immersed minimal annulus, with two planar boundary
curves along which the surface meets these planes with constant contact
angle, is part of the catenoid.
\end{abstract}
\pagestyle{plain}
\bigskip
\section{Introduction}
\bigskip
$\quad$ The catenoid is the simplest minimal surface in
$\mathbb{R}^3$ except the plane. It is obtained by revolving the
catenary about an axis. The catenoid has been characterized by
many authors. For instance:
 \begin{itemize}
\item[(1)] The catenoid is the only nonplanar minimal surface
which is a surface of revolution (Bonnet \cite{nit}). \item[(2)]
The catenoid is the only complete embedded minimal surface of
total curvature $-4\pi$ (Osserman \cite{oss}). \item[(3)] A
complete minimal surface with two annular ends and of finite total
curvature is the catenoid (Schoen \cite{sch}). \item[(4)] A
complete embedded nonplanar minimal surface of finite total
curvature and genus zero is the catenoid (L\'{o}pez and Ros
\cite{lr}).
\end{itemize}
\indent For more interesting characterizations see
\cite{ch}, \cite{coll} and \cite{n}.\\
\indent On the other hand, one can also characterize a proper subset of the catenoid.\\
 In 1869, Enneper \cite{nit} proved that a compact nonplanar minimal surface which is generated
 by one-parameter family of circles is part of the catenoid or Riemann's
 example. In 1956,  Shiffman \cite{max} proved that a minimal annulus bounded by two
horizontal circles is foliated by horizontal circles.

In this paper we will also characterize a proper subset of the
catenoid. Our characterizations involves the hypothesis of a
constant contact angle along the boundary of the minimal surface as
follows:

\begin{thm} Let $\Sigma$ be
an immersed minimal annulus such that $\partial\Sigma$ consists of
two $C^{2,\alpha}$ planar Jordan curves $\Gamma_1$ and $\Gamma_2$.
If $\Sigma$ makes a constant contact angle with a plane $\Pi_i$
along $\Gamma_i$, $i=1,2$, $\Pi_1\neq\Pi_2$, then $\Sigma$ is part
of the catenoid.
\end{thm}

\begin{thm} Let $\Sigma$ be an immersed minimal surface with boundary
 and let $\Gamma$ be one component of
$\partial\Sigma$. If $\Gamma$ is a circle and $\Sigma$ meets a
plane along $\Gamma$ at a constant angle, then $\Sigma$ is part of
the catenoid.
\end{thm}

Note that in Theorem 1.1 it is not necessary to assume that the
planes $\Pi_1$ and $\Pi_2$ are parallel.  Also it should be
mentioned that Wente \cite{we1} proved every embedded annular
capillary surface in a slab is a surface of revolution. But he also
constructed many examples of immersed non-zero constant mean
curvature (henceforth abbreviated as CMC) annular capillary surfaces
lying in a slab which are not surfaces of revolution \cite{we2}.

The author would like to express his gratitude to professor 
J. Choe for his guidance and encouragement.
\section{Preliminaries}
\bigskip
$\quad$ First, we review the Hopf differential. Let $\Sigma$ be an
annulus in $\mathbb{R}^3$ which is the image of a conformal
immersion $X$ of a planar annulus $ A = \{(u,v) \in {\mathbb{R}^2} :
{1}/{R}\leq u^2 +v^2 \leq R, R>1 \}$. Suppose $u$ and $v$ are the
isothermal coordinates on $A$ determined by $X$. We can write the
first fundamental form and the second fundamental form of $\Sigma$
as follows
\begin{equation*}\label{eq:g}
 I_{X} = E({du}^2 + {dv}^2),
\end{equation*}
\begin{equation*}\label{eq:yu}
{ {II}}_{ X} = L{du}^2 + 2M{du}{dv} + N{dv}^2.
\end{equation*}
 The {Hopf differential} is the quadratic differential defined by $\Phi dw^2$, $\Phi = ({L - N})/{2}
 -i{M}$, $w=u+i v$. Then the Codazzi equation implies the following lemma.
\begin{lem}(See \cite{c} or \cite{hopf}.)
$\Phi$ is holomorphic on a CMC surface.
\end{lem}
Second, we review some properties of umbilic points of a CMC
surface. Umbilic points of $\Sigma$ are the zeros of $\Phi$. Lines
of curvature of $\Sigma$ flow smoothly except at umbilic points.
They rotate sharply around at an umbilic point. So we can define the
rotation index of the lines of
curvature at interior umbilic points.\\
\indent Now we extend the rotation index to a boundary point. Let
$p\in\partial\Sigma$ be a boundary point. We choose $X : D^{+}
\rightarrow \Sigma $ which is a conformal immersion of a half disk $
D^{+} = \{(u,v)\in D : u^2+v^2\leq1, v\geq 0 \}$ into the regular
surface $\Sigma$ mapping the diameter $l$ of $D^{+}$ into
$\partial\Sigma$. The lines of curvature of $\Sigma$ can be pulled
back by $X$ to a line field on $ D^{+}$. If $X(l)$ is a line of
curvature of $\Sigma$, then this line field can be extended smoothly
to a line field $F$ on $ D$ by reflection about the diameter $l$. It
is clear that $F$ has the well-defined rotation index at $X^{-1}(p)$
and furthermore, the rotation index does not depend on the choice of
immersion $X$. So we can naturally define the rotation index of the
lines of curvature at the umbilic point $p\in \partial\Sigma$ to be
half the rotation index of $F$ at $X^{-1}(p)$.
\begin{lem} (\cite{c}, Lemma 2)\\
Let $\Sigma$ be a non-totally umbilic immersed CMC surface which is
of class $C^{2,\alpha}$ up to and including the boundary
$\partial\Sigma$. If the $\partial\Sigma$ are lines of curvature,
then we have the following.\\
(a) The boundary umbilic points of $\Sigma$ are isolated.\\
(b) At an interior umbilic point the rotation index of lines of
curvature is not bigger than $-{1}/{2}$. \\
(c) At a boundary umbilic point the rotation index of lines of
curvature is not bigger than $-{1}/{4}$.
\end{lem}
We now recall Bj\"orling's theorem (see \cite{di}). Let
$c:[a,b]\rightarrow \mathbb{R}^3$ be any real analytic curve and
$n:[a,b]\rightarrow \mathbb{S}^2$ be any real analytic vector field
perpendicular to the tangent vector of the curve $c(t)$.  By the
analyticity of $c$, there are unique analytic extensions
$c:[a,b]\times(-\varepsilon,\varepsilon)\rightarrow \mathbb{C}^3$,
and $n:[a,b]\times(-\varepsilon,\varepsilon)\rightarrow
\mathbb{C}^3$, where $\varepsilon$ is a small enough positive
number. Using these extensions, we define the unique immersion of surface as
follows
\begin{eqnarray} \label{eq:c1}
X(z)=\re \Big( c(z)-i\int^{z}_{0}n(w)\times c'(w) dw\Big),
\end{eqnarray}
where $z \in [a,b]\times(-\varepsilon,\varepsilon).$ By a
straightforward computation, this is a minimal immersion that
extends $c$ and $n$ in the sense that for $t \in [a,b]\times \{0\}$,
$X(t)=c(t)$ and $n(t)$ is the surface normal.
\section{Proof of the theorems}
\bigskip
\textit{Proof of theorem 1.1.}
\textbf{Step 1.} We claim that $\Pi_1$ and $\Pi_2$ are parallel.\\
 Let $X  :  A = \{(u,v) \in {\mathbb{R}^ 2} : {1}/{R}\leq u^2
 +v^2 \leq R \} \rightarrow \mathbb{R}^3$ be a conformal immersion of
$\Sigma$. Because $\Sigma$ is a minimal surface, $\triangle X = 0 $
on $A$, where $\triangle$ is the Laplace-Beltrami operator on
$\Sigma$. So we have
\begin{equation}\label{eq:rt}
 0 = \int_{\Sigma}^{} \triangle X  d A
=\int_{\Gamma_1}^{}\nu_1 ds + \int_{\Gamma_2}^{}\nu_2 ds,
\end{equation}
where $\nu_i$ denotes the outward pointing unit conormal vector
along  $\Gamma_i$, $i=1,2$. By the boundary maximum principle
\cite{hopf}, $\theta_1\neq0,\pi$ (see, Figure 1). Let $\Omega_1\subset\Pi_1$ be the
domain bounded by $\Gamma_1$. The projection of $\nu_1$ to $\Pi_1$
is a normal vector field of $\Gamma_1$, and it is constant length.
By the divergence theorem on $\Omega_1$, we see that the non-zero
vector $\int_{\Gamma_1}^{}\nu_1 ds$ is perpendicular to $\Pi_1$.
Similarly, the non-zero vector $\int_{\Gamma_2}^{}\nu_2 ds$ is also
perpendicular to $\Pi_2$. (\ref{eq:rt}) implies that the two vectors
are linearly dependent.
So $\Pi_1$ and $\Pi_2$ are parallel.\\
\indent\textbf{Step 2.} We claim that both $\Gamma_{1}$ and
$\Gamma_{2}$ are convex.\\
The Terquem-Joachimsthal theorem \cite{spi} says that if
$\Gamma=\Sigma_1 \cap \Sigma_2$ is a line of curvature in
$\Sigma_1$, then $\Gamma$ is also a line of curvature in $\Sigma_2$
if and only if $\Sigma_1$ and $\Sigma_2$ intersect at a constant
angle along $\Gamma$. Since $\Gamma_1$ is a line of curvature of
$\Pi_1$ and $\Sigma$ meets $\Pi_1$ in a constant contact angle along
the $\Gamma_1$, $\Gamma_1$ is also a line of curvature of $\Sigma$.
Similarly, $\Gamma_2$ is a line of curvature of $\Sigma$.

Let $\kappa_1$, $\kappa_2$ be the principal curvature of $\Sigma$
along $\Gamma_1$, $\Gamma_2$ respectively. We want to show that neither
$\kappa_1$ nor $\kappa_2$ has zeros.

First, let us suppose that $\kappa_{i}$ has zeros at finite points
$p_{{i}_{j}},~j=1,...,m_{i}$ on $\Gamma_{i}$, $i=1,2.$ Then
$p_{{i}_{j}},~j=1,...,m_{i}$, $i=1,2$ are the boundary umbilic
points. Let $q_{k},~k=1,...,n$ be the interior umbilic points. By
the Poincar\'{e}-Hopf theorem and Lemma 2.2, we have
\begin{equation*}\label{eq:ty}
  \chi(\Sigma)= 0 = \sum_{p=p_{i_{j}},q_{k}}I(p)\leq
 \sum_{i}\sum_{{j}}(-\frac{1}{4}) +\sum_{k}(-\frac{1}{2}) < 0,
\end{equation*}
where $\chi(\Sigma)$ is the Euler characteristic of $\Sigma$ and
$I(p)$ is the rotation index at $p$. Therefore neither $\kappa_{1}$ nor
$\kappa_{2}$  has zeros.

Second, suppose either  $\kappa_{1}$ or $\kappa_{2}$ has
zeros at an infinite number of points. By lemma 2.1, $\Phi$ is a holomorphic
function. Since $\Gamma_1,$ $\Gamma_2$ are compact sets, the Hopf
differential $\Phi$ of $\Sigma$ is identically zero on the one or
both of $\Gamma_1$ and $\Gamma_2$. So $\Phi$ is identically zero.
This means that $\Sigma$ is a planar annulus. Since
$\Pi_1\neq\Pi_2$, this case cannot happen.

Hence both $\kappa_1$ and $\kappa_2$ cannot have zeros.

Let $\widetilde{\kappa}_i$ be the curvature  of $\Gamma_i$ in
$\Pi_i$. Since $\kappa_{i}=\widetilde{\kappa}_i\sin\theta_i$,
$\widetilde{\kappa}_i$ cannot be zero.  So both
$\Gamma_1$ and $\Gamma_2$ are convex. In fact they are strictly convex.\\
\indent \textbf{Step 3.} We claim that $\Sigma$ is  part of the catenoid.\\
By Step 2 the total curvature of $\partial\Sigma$ is $4\pi$. So
by \cite{w}, $\Sigma$ is embedded (see, Figure 1).\\
Hence we know that $\Sigma$ is a surface of revolution by \cite{we1}.
But we give a sketch of the proof for the sake of completeness.\\
\indent Let $v$ be a unit vector in $\Pi_1$.
 Now we apply the Alexandrov reflection principle \cite{hopf} with the one-parameter family of plane $\Pi_{v,t}$ which is  orthogonal to $v$. Increasing $t$ one gets a first plane $\Pi_{v,\overline{t}}$ that reach $\Sigma$; that is $\Pi_{v,\overline{t}}\cap \Sigma\neq\emptyset$, but if $t<\overline{t}$ then $\Pi_{v,t}\cap \Sigma=\emptyset$. Let $\Sigma_t$ be the part of $\Sigma$ lying in $\Pi_{v,s}$, $s<t$ and let $\Sigma_{\Pi_{v,t}}$ be the symmetry of $\Sigma_t$ about $\Pi_{v,t}$. First, let us assume the first  touching between $\Sigma$ and the reflected surface
 $\Sigma_{\Pi_{v,t_0}}$ by $\Pi_{v,t_0}$ for some $t_0$ occurs at an interior
 point. By the Hopf maximum principle, $\Sigma$ is symmetric with respect to $\Pi_{v,t_0}$.
 Second, let us assume the first touching between $\Sigma$ and
 $\Sigma_{\Pi_{v,t_1}}$ for some $t_1$ occurs at a boundary point. Since the contact angle is constant,
 the normal vector of  $\Sigma$ at the touching point coincides with that of the reflected surface
 $\Sigma_{\Pi_{v,t_1}}$ at the touching point. So we can use the
 Hopf boundary maximum principle, then $\Sigma$ is symmetric with respect to
 $\Pi_{v,t_1}$. Otherwise, let us assume the first touching between $\Sigma$ and
 $\Sigma_{\Pi_{v,t_2}}$ for some $t_2$ occurs at a corner. Similarly the second case, the normal vector of
 $\Sigma$ at the corner coincides with that of the reflected surface $\Sigma_{\Pi_{v,t_2}}$ at the corner.
 Because of constant of the contact angle we can apply the Serrin's boundary point lemma at a corner (Lemma 2.6 of \cite{we1} or \cite{se}).
 Then $\Sigma$ is symmetric with respect to  $\Pi_{v,t_2}$.\\
 Since the first touching must occur at an interior point, a boundary point or a corner, $\Sigma$ is symmetric with respect to $\Pi_{v,T}$.\\
 \indent By using the reflection principle for another unit
 vector $w$ in $\Pi_1$, we get that $\Sigma$ is symmetric with respect to $\Pi_{w,T_{1}}$. Hence, we conclude that $\Sigma$ is a rotational surface, $i.e.$, it is part of the catenoid.  $\Box$\\

\begin{figure}
\begin{center}
\includegraphics[height=6cm, width=11cm]{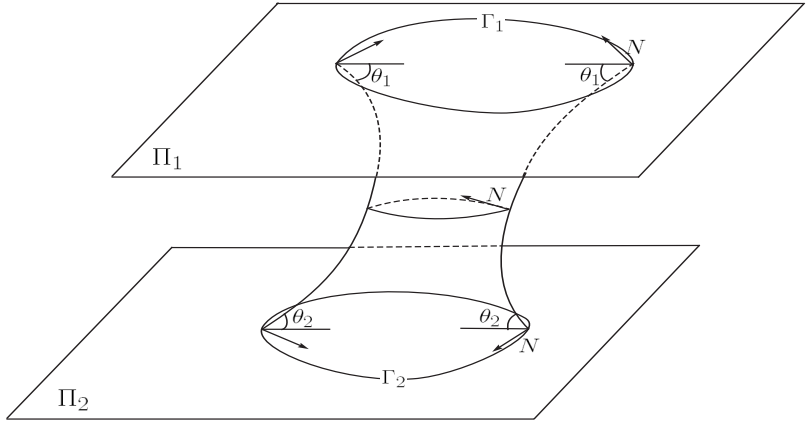}
\end{center}
\caption{Minimal annulus with constant contact angles.}
\end{figure}

\textit{Proof of theorem 1.2.} Without loss of generality denote the
circle by $c(t)=(\cos t,\sin t,0)$. Since the surface has constant
contact angle $\theta$ along the $c$, the surface normal vector
becomes $n(t)=(\sin \theta \cos t,\sin\theta\sin t,\cos\theta)$.
Then Bj\"orling's formula (\ref{eq:c1}) yields the unique minimal surface
$$X(z)=(X_{1}(z),X_{2}(z),X_{3}(z))=\re\Big((\cos z,\sin z,0)-i\int^{z}_{0}\Big(-\cos\theta \cos w, \cos\theta \sin w, \sin\theta\Big)dw\Big)$$
$$=(\cos u \cosh v-\cos \theta\sin u \sinh v , \sin u \cosh v + \cos \theta \cos u \sinh v , \sin\theta v).$$
For each $X_{3}$, $X^2_{1}+X^2_{2}=\constant$. This means that the
minimal surface is foliated by coaxial circles. So $\Sigma$ is part
of the catenoid. $\Box$

\section{Remarks}
\indent Lemma 2.1, Lemma 2.2 and the Poincar\'{e}-Hopf theorem hold
for CMC surfaces. So we drive a characterization of a sphere after an
additional assumption.
\begin{cor} Let $\Sigma$ be
an immersed non-zero CMC annulus such that $\partial\Sigma$ consists
of two $C^{2,\alpha}$ planar Jordan curves $\Gamma_1$ and
$\Gamma_2$. If $\Sigma$ makes a constant contact angle with a plane
$\Pi_i$ along $\Gamma_i$, $i=1,2$, $\Pi_1\neq\Pi_2$. In addition,
$\Sigma$ has at least one umbilic point. Then $\Sigma$ is part of a
sphere.
\end{cor}
\begin{proof}
If $\Sigma$ has only finite umbilic points then it is contradiction
to the Poincar\'{e}-Hopf theorem.  So  $\Sigma$ has infinitely many
umbilic points.  Since the Hopf differential $\Phi$ is holomorphic,
it is identically zero. Hence, $\Sigma$ is part of a
sphere.
\end{proof}

So far, we have considered minimal annuli with boundary curves lying in a pair of planes. In
case of minimal annuli with boundary curves lying on a sphere, there is a well-known open
problem.
\begin{p}(\cite{sou})
Let $\Gamma_1$, $\Gamma_2$ be  ${C}^{2,\alpha}$  Jordan curves on a
sphere. Show that if $\Sigma$ is an immersed minimal annulus meeting
constant contact angles with the sphere along the boundary
$\Gamma_1$, $\Gamma_2$, then $\Sigma$ is part of the catenoid.
\end{p}
\begin{rem}
By Theorem 1.2, Problem 4.2 is true if one of the boundary curves is
a circle.
\end{rem}

\vspace{1cm} \noindent Department of Mathematics, Seoul National
University,
Seoul 151-742, Korea\\
\tt{E-mail: jcpyo@snu.ac.kr}\bigskip

\end{document}